\documentclass[a4paper,english,12pt]{article}
\usepackage[logonly]{trace}
\usepackage{babel}
\usepackage{amsfonts, amsmath, amssymb}
\usepackage{graphicx}
\usepackage{pstricks}
\usepackage{psfig}
\usepackage{multirow}
\usepackage{fancyhdr}
\usepackage{vmargin, fancybox}

\newcommand{\mathsym}[1]{{}}

\def\N{\textrm{I\kern-0.21emN}}

\def\R{\textrm{I\kern-0.21emR}}
\def\Q{\textrm{l\kern-0.5emQ}}

\begin{document}

\title{Structural stability of finite dispersion-relation preserving schemes}

\author{Claire David\footnotemark[2] \footnote{Corresponding author:
david@lmm.jussieu.fr; fax number: (+33) 1.44.27.52.59.} ,  Pierre
Sagaut\footnotemark[2] \\ \small{Universit\'e Pierre et Marie Curie-Paris 6}  \\
\footnotemark[2] \small{ Institut Jean Le Rond d'Alembert, UMR CNRS 7190,} \\
\small{Bo\^ite courrier $n^0 162$, 4 place Jussieu, 75252 Paris,
cedex 05, France}}


\maketitle

\begin{abstract}
The goal of this work is to determine classes of travelling
solitary wave solutions for a differential approximation of a
finite difference scheme by means of a hyperbolic ansatz. It is
shown that spurious solitary waves can occur in finite-difference
solutions of nonlinear wave equation. The occurance of such a
spurious solitary wave, which exhibits a very long life time,
results in a non-vanishing numerical error for arbitrary time in
unbounded numerical domain. Such a behavior is referred here to
has a structural instability of the scheme, since the space of
solutions spanned by the numerical scheme encompasses types of
solutions (solitary waves in the present case) that are not
solution of the original continuous equations. This paper extends
our previous work about classical schemes to dispersion-relation
preserving schemes.
\end{abstract}

\section{Introduction: The DRP scheme}
\label{DRP}

\noindent The Burgers equation:

\begin{equation}
\label{Burgers}  u_t + c\, u\, u_x - \mu \,u_{xx} = 0,
\end{equation}

\noindent $c$, $\mu$ being real constants, plays a crucial role in
the history of wave equations. It was named after its use by
Burgers \cite{burger1} for studying
turbulence in 1939.\\

\noindent $i$, $n$ denoting natural integers, a linear finite
difference scheme for this equation can be written under the form:
\begin{equation} \label{scheme} \displaystyle \sum \alpha_{lm}\,u_{l}^{m}=0
              \end{equation}

\noindent where:
\begin{equation}
{u_l}^m=u\,(l\,h, m\,\tau)
\end{equation}
\noindent  $l\, \in \, \{i-1,\, i, \, i+1\}$, $m \, \in \,
\{n-1,\, n, \, n+1\}$, $j=0, \, ..., \, n_x$, $n=0, \, ..., \,
n_t$. The $ \alpha_{lm}$ are real coefficients, which depend on the mesh size $h$, and the time step $\tau$.\\
The Courant-Friedrichs-Lewy number ($cfl$) is defined as $\sigma = c \,\tau / h$ .\\
\noindent A numerical scheme is  specified by selecting
appropriate values of the coefficients $ \alpha_{lm}$. Then,
depending on them, one can obtain optimum schemes, for which the
error will be minimal.\\

 \noindent $m$ being a strictly positive integer, the first derivative $\frac{\partial
u}{\partial x}$ is approximated at the $l^{th}$ node of the
spatial mesh by:

\begin{equation}\label{approx}
 (\, \frac{\partial u}{\partial
x}\,)_l  \simeq
   \displaystyle \sum_{k=-m}^m \gamma_{k}\,u_{i+k}^n
\end{equation}
\noindent Following the method exposed by C. Tam and J. Webb in
\cite{Tam}, the coefficients $\\gamma_{k}$ are determined
requiring the Fourier Transform of the finite difference scheme
(\ref{approx}) to be a close approximation of the partial
derivative $ (\, \frac{\partial u}{\partial x}\,
)_l$.\\
\noindent (\ref{approx}) is a special case of:

\begin{equation}\label{approx_Cont}
 (\, \frac{\partial u}{\partial
x}\,)_l  \simeq \displaystyle \sum_{k=-m}^m \gamma_{k}\,u(x+k\,h)
\end{equation}

\noindent where $x$ is a continuous variable, and can be recovered
setting $x=l\,h$.\\
\noindent Denote by $\omega$ the phase. Applying the Fourier
transform, referred to by $\,\widehat{\, }$ , to both sides of
(\ref{approx_Cont}), yields:

\begin{equation}
\label{Wavenb}
 j\, \omega \, \widehat{u}  \simeq \displaystyle \sum_{k=-m}^m \gamma_{k}\,e^{\,j\,k\,\omega\,h}\, \widehat{u}
\end{equation}
\noindent  $j$ denoting the complex square root of $-1$.\\




\noindent Comparing the two sides of (\ref{Wavenb}) enables us to
identify the wavenumber $ \overline{\lambda}$ of the finite
difference scheme (\ref{approx}) and the quantity
$\frac{1}{j}\,{\displaystyle \sum_{k=-m}^m
\gamma_{k}\,e^{\,j\,k\,\omega\,h}}$, i. e.: \noindent The
wavenumber of the finite difference scheme (\ref{approx}) is thus:

\begin{equation}
 \overline{\lambda}=-\,j\, \displaystyle \sum_{k=-m}^m \gamma_{k}\,e^{\,j\,k\,\omega\,h}
\end{equation}

\noindent To ensure that the Fourier transform of the finite
difference scheme is a good approximation of the partial
derivative $ (\, \frac{\partial u}{\partial x}\, )_l$ over the
range of waves with wavelength longer than $4\,h$, the a priori
unknowns coefficients $\gamma_{k}$ must be choosen so as to
minimize the integrated error:

\footnotesize
\begin{equation}\begin{array}{rcl}
 {\mathcal E} &=&\int_{-\frac{\pi}{2}}^{\frac{\pi}{2}} | \lambda \,h- \overline{\lambda}
 \,h|^2\,d(\lambda \,h)\\
 &=&\int_{-\frac{\pi}{2}}^{\frac{\pi}{2}} | \lambda \,h+j\, \displaystyle \sum_{k=-m}^m \gamma_{k}\,e^{\,j\,k\,\omega\,h} \,h|^2\,d(\lambda \,h)\\
 &=& \int_{-\frac{\pi}{2}}^{\frac{\pi}{2}} | \zeta+j\, \displaystyle \sum_{k=-m}^m \gamma_{k}\,\left  \lbrace \cos (\,k\,\zeta)+j\,\sin(\,k\,\zeta) \right \rbrace \,
 |^2\,d \zeta \\
 &=& \int_{-\frac{\pi}{2}}^{\frac{\pi}{2}} \left \lbrace \left [ \zeta- \displaystyle \sum_{k=-m}^m \gamma_{k} \,\sin(\,k\,\zeta) \right ]^2+
  \left [   \displaystyle \sum_{k=-m}^m \gamma_{k}\, \cos (\,k\,\zeta)  \right ]^2\, \right \rbrace  \,d
  \zeta \\
   &=& 2 \, \int_0^{\frac{\pi}{2}} \left \lbrace \left [ \zeta- \displaystyle \sum_{k=-m}^m \gamma_{k} \,\sin(\,k\,\zeta) \right ]^2+
  \left [   \displaystyle \sum_{k=-m}^m \gamma_{k}\, \cos (\,k\,\zeta)  \right ]^2\, \right \rbrace  \,d
  \zeta \\
      \end{array}
\end{equation}

\normalsize

\noindent The conditions that ${\mathcal E}$ is a minimum are:

\begin{equation}
 \frac {\partial{\mathcal E}}{\partial \gamma_i} =0 \,\,\,  , \,\,\,
i=-m, \ldots , \, m
\end{equation}

\noindent i. e.:

\begin{equation}
\label{RelDer}
  \int_0^{\frac{\pi}{2}} \left \lbrace
 \,-\,\zeta\,\sin(\,i\,\zeta)\, + \displaystyle \sum_{k=-m}^m \gamma_{k}
\,\cos\left(\,(k-i)\,\zeta \right)    \right \rbrace  \,d
  \zeta =0
\end{equation}

\noindent Changing $i$ into $-i$, and $k$ into $-k$ in the
summation yields:

\begin{equation}
  \int_0^{\frac{\pi}{2}} \left \lbrace
 \,\,\zeta\,\sin(\,i\,\zeta)\, + \displaystyle \sum_{k=-m}^m \gamma_{-k}
\,\cos\left(\,(-k+i)\,\zeta \right)    \right \rbrace  \,d
  \zeta =0
\end{equation}

\noindent i. e.:

\begin{equation}
  \int_0^{\frac{\pi}{2}} \left \lbrace
 \,\,\zeta\,\sin(\,i\,\zeta)\, + \displaystyle \sum_{k=-m}^m \gamma_{-k}
\,\cos\left(\,(k-i)\,\zeta \right)    \right \rbrace  \,d
  \zeta =0
\end{equation}

\noindent Thus:

\begin{equation}
\label{Int1}
  \int_0^{\frac{\pi}{2}}  \displaystyle \sum_{k=-m}^m \left \lbrace
  \gamma_{-k} + \gamma_{k}  \right \rbrace
\,\cos\left(\,(k-i)\,\zeta \right)    \,d
  \zeta =0
\end{equation}

\noindent which yields:

\begin{equation}
 \frac{\pi}{2} \,\left \lbrace
  \gamma_{-i} + \gamma_{i}  \right \rbrace +\displaystyle \sum_{k \neq  i,\, k=-m}^m
\left \lbrace \frac {
  \gamma_{-k} + \gamma_{k} }{k-i} \right \rbrace
\,\sin\left(\,(k-i)\,\frac{\pi}{2} \right)    \  =0
\end{equation}

\noindent which can be considered as a linear system of $2\,m+1$
equations, the unknowns of which are the $ \gamma_{-i} +
\gamma_{i} $, $i=-m,\,\ldots, \,m$. The determinant of this system
is not equal to zero, while it is the case of its second member:
the Cramer formulae give then, for $i=-m,\,\ldots, \,m$:

\begin{equation}
  \gamma_{-i} + \gamma_{i}    =0
\end{equation}

 \noindent or:

\begin{equation}
\label{RelCoeffDRP1}
  \gamma_{-i} =- \gamma_{i}
\end{equation}

 \noindent For $i = 0$, one of course obtains:

\begin{equation}
   \gamma_0=0
\end{equation}

 \noindent All this ensures:

\begin{equation}
\label{RelCoeffDRP} \displaystyle \sum _{k=-m}^m \gamma_{k}  =0
\end{equation}

 \noindent The values of the $\gamma_k$ coefficients are obtained by substituting relations (\ref{RelCoeffDRP1}) into (\ref{RelDer}):

\begin{equation}
\label{RelCoeffDRP} \displaystyle \sum _{k=-m}^m \gamma_{k}  =0
\end{equation}

\noindent $m$ being a strictly positive integer, a ${2m+1}$-points
\textit{DRP} scheme (\cite{Tam}) is thus given by:

\footnotesize

\begin{equation}
\label{DRP}
\begin{aligned}
&-u_{i}^{n+1}+u_{i}^{n } +\frac{ \tau }{h}\, \displaystyle \sum
_{k=-m}^m \gamma_{k}\,u_{i+k}^{n }  =0
   \end{aligned}
\end{equation}

\normalsize

 \noindent where the $\gamma_{k}$, $k\in \{-m,m\}$ are the coefficients of the considered
 scheme, and satisfy the relations (\ref{RelCoeffDRP1}).

\noindent Considering again the ${u_l}^m$ terms as functions of
the mesh size $h$ and time step $\tau$, expanding them at a given
order by means of their Taylor series expansion, and neglecting
the $o(\tau^p)$ and $o({h}^q)$ terms, for given values of the
integers $p$, $q$, leads to the following differential
approximation (see \cite{Shokin}):

\footnotesize

\begin{equation}
\label{Approx}
\begin{aligned}
&-u_{i}^{n+1}+u_{i}^{n } +\frac{ \tau }{h}\,\displaystyle \sum
_{k=-m}^m \gamma_{k}\, \mathcal {F}_{i+k}^{n } \left (u,\frac
{\partial^r u}{\partial x^r},\frac {\partial^s u}{\partial
t^s},h,\tau  \right ) =0
   \end{aligned}
\end{equation}

\normalsize

\noindent where $\mathcal {F}_{i+k}^{n }$ denotes the function of
$u$, $\frac {\partial^r
u}{\partial x^r},\frac {\partial^s u}{\partial t^s}$, $h$, $\tau$ obtained by means of the above Taylor expansion, $r$, $s$ being integers.\\

\noindent For sake of simplicity, a non-dimensional form of
Eq. (\ref{Approx}) will be used:

\begin{equation}
\label{ApproxAdim}
\begin{aligned}
&-\tilde{u}_{i}^{n+1}+\tilde{u}_{i}^{n } +   \displaystyle \sum
_{k=-m}^m \gamma_{k}\, \mathcal {F}_{i+k}^{n }
\left(\tilde{u},\frac {\partial^r \tilde{u}}{\partial
x^r},\frac{\partial^s u}{\partial \tilde{t}^s} \right)  =0
   \end{aligned}
\end{equation}

\noindent Depending on this differential approximation
(\ref{ApproxAdim}), solutions,
as solitary waves, may arise.\\

\noindent The paper is organized as follows. DRP schemes are
analyzed in section \ref{Analyse}. The general method is exposed
in Section \ref{Solitons}. Classical DRP schemes are studied in
section \ref{Solutions}, where it is shown that out of the two
studied schemes, only one leads to solitary waves. A related class
of travelling wave solutions of equation (\ref{Approx}) is thus
presented, by using a hyperbolic ansatz.

\section{Analysis of \textit{DRP} schemes}
\label{Analyse}

\noindent Consider ${u_{i}}^{n }$ as a function of the time step
$\tau$, and expand it at the second order by means of its Taylor
series:

\begin{equation}
{u_i}^{n+1}=u\,(i\,h, (n+1)\,\tau)=u\,(i\,h,
n\,\tau)+\tau\,u_t\,(i\,h,
n\,\tau)+\frac{{\tau}^2}{2}\,u_{tt}\,(i\,h, n\,\tau)+o({\tau}^2)
\end{equation}

\noindent It ensures:

\begin{equation}
\label{Prop1} \frac{{u_i}^{n+1}-{u_i}^{n}}{\tau}=u_t\,(i\,h,
n\,\tau)+\frac{{\tau}}{2}\,u_{tt}\,(i\,h, n\,\tau)+o({\tau})
\end{equation}

\noindent In the same way, for $k \in \{-m,m\}$, consider
${u_{i+k}}^{n }$ as a function of the mesh size $h$, and expand it
at the fourth order by means of its Taylor series expansion:

\begin{equation}
\label{Prop_a}  \scriptsize  \begin{aligned} {u_{i+k}}^{n}  = &
u\,((i+k)\,h, n\,\tau)\\& =  u\,({ }\,h,
n\,\tau)+k\,h\,u_x\,(i\,h, n\,\tau)\\&+\frac{ k^2\,h
^2}{2}\,u_{xx}\,( i\,h, n\,\tau)+\frac{ k^4\,h
^4}{4!}\,u_{xxxx}\,(i\,h, n\,\tau)+o({h}^4)
\end{aligned}
\end{equation}

\bigskip \noindent  Equation (\ref{Approx}) can thus be written as:

\normalsize

\footnotesize

\begin{equation}
\label{BurgersEq}   \scriptsize \begin{aligned} & -u_t\,(i\,h,
n\,\tau)-\frac{{\tau}}{2}\,u_{tt}\,(i\,h,
n\,\tau)+o({\tau})\\&+\frac{ \tau }{h}\,\displaystyle \sum
_{k=-m}^m \gamma_{k}\, \left \lbrace u\,(i\,h,
n\,\tau)+k\,h\,u_x\,(i\,h, n\,\tau) +\frac{ k^2\,h
^2}{2}\,u_{xx}\,( i\,h, n\,\tau)+\frac{{ k}^4\,h
^4}{4!}\,u_{xxxx}\,(i\,h, n\,\tau)+o({h}^4) \right \rbrace =0
\end{aligned}
\end{equation}

\normalsize

\noindent i. e., at $x=i\,h$ and $t=n \,\tau$:

\footnotesize

\begin{equation}
\scriptsize  \begin{aligned} & -u_t -\frac{{\tau}}{2}\,u_{tt}
+o({\tau}) +\frac{ \tau }{h}\,\displaystyle \sum _{k=-m}^m
\gamma_{k}\, \left \lbrace u +k\,h\,u_x  +\frac{ k^2\,h
^2}{2}+u_{xx} \frac{{ k}^4\,h ^4}{4!}\,u_{xxxx} +o({h}^4) \right
\rbrace =0
\end{aligned}
\end{equation}

\normalsize

\noindent (\ref{RelCoeffDRP}) ensures then:

\footnotesize

\begin{equation}
\scriptsize  \begin{aligned} & -u_t -\frac{{\tau}}{2}\,u_{tt}
+o({\tau}) +\frac{ \tau }{h}\,\displaystyle \sum _{k=-m}^m
k\,\gamma_{k} \, \left \lbrace  h\,u_x  +o({h}^4) \right \rbrace
=0
\end{aligned}
\end{equation}

\normalsize

\noindent The related first differential approximation of the
Burgers equation (\ref{Burgers}) is thus obtained neglecting the
$o(\tau)$ and $o({h}^2)$ terms, yielding:

\begin{equation} \label{Eq}
\scriptsize  \begin{aligned} & -u_t -\frac{{\tau}}{2}\,u_{tt}
  + { \tau } \,\displaystyle \sum _{k=-m}^m
k\,\gamma_{k} \,u_x   =0
\end{aligned}
\end{equation}

\noindent For sake of simplicity, this latter equation can be
adimensionalized in the following way:\\
\noindent set:

\begin{equation}
\label{Chgt_Var_Adim} \left \lbrace
\begin{array}{rcl}
u &  =   & U_0\,\tilde {u}\\
t &  =   & \tau_0\,\tilde {t}\\
x &  =   & h_0\,\tilde {x}
\end{array}
\right. \end{equation}

\noindent where:

\begin{equation}
\label{U0} U_0=\frac {h_0}{\tau_0}\end{equation}

\noindent In the following,  $ Re_h$ will denotes the mesh Reynolds number, defined as:

\begin{equation}
\label{Re}
 Re_h=\frac{U_0\,h}{\mu}
\end{equation}

\noindent For $j \in \N$, the change of variables
(\ref{Chgt_Var_Adim}) leads to:

\begin{equation}
\label{Chgt_Var} \left \lbrace
\begin{array}{rcl}
u_t &  =   & \frac {U_0}{\tau_0} \,\tilde {u}_{\tilde {t}}\\
u_{t^j} &  =   & \frac {U_0}{\tau_0^j} \,\tilde {u}_{{\tilde
{t}}^j}\\
 u_{x^j} &  =   & \frac {U_0}{h_0^j} \,\tilde {u}_{{\tilde
{x}}^j }
\end{array}
\right. \end{equation}

\noindent (\ref{Eq}) becomes:

\begin{equation}
\scriptsize  \label{Eq2} \begin{aligned} &  -\frac {U_0}{\tau_0}
\,\tilde {u}_{\tilde {t}} -\frac{{\tau}}{2}\,\frac {U_0}{\tau_0^2}
\,\tilde {u}_{{\tilde {t}}{\tilde {t}}}
  + 2\,{ \tau } \,\displaystyle \sum _{k=1 }^m
k\,\gamma_{k}  \,\frac {U_0}{h_0 } \,\tilde {u}_{{\tilde {x}}  }
=0
\end{aligned}
\end{equation}

\noindent Multiplying (\ref{Eq2}) by $\frac {\tau_0}{U_0}$ yields:

\begin{equation}
\label{Eq3} \scriptsize  \begin{aligned} &  - \tilde {u}_{\tilde
{t}} -\frac{\tau}{2\,\tau_0 } \,\tilde {u}_{{\tilde {t}}{\tilde
{t}}}
  +2\,\frac{ \tau\, {\tau_0} }{h_0}\,\displaystyle \sum _{k=1}^m
k\,\gamma_{k} \,h\,  \tilde {u}_{{\tilde {x}} }      =0
\end{aligned}
\end{equation}

\noindent For $h=h_0$, due to $\sigma=\frac {U_0\,\tau}{h}$, Eq.
(\ref{Eq3}) becomes:

\begin{equation}
\scriptsize  \label{EqEq} \begin{aligned} &  - \tilde {u}_{\tilde
{t}} -\frac{\tau}{2\,\tau_0 } \,\tilde {u}_{{\tilde {t}}{\tilde
{t}}}
  +2\,\frac{ \tau\,h_0  }{\mu \,Re_h }\,\displaystyle \sum _{k=1}^m
k\,\gamma_{k} \,     \,\tilde {u}_{{\tilde {x}}  }   =0
\end{aligned}
\end{equation}

\noindent which simplifies in:

\begin{equation}
\scriptsize  \label{EqEq} \begin{aligned} &  - \tilde {u}_{\tilde
{t}} -\frac{\sigma}{2  } \,\tilde {u}_{{\tilde {t}}{\tilde {t}}}
  +\frac{ 2\,\sigma  }{\mu\,Re_h} \,\displaystyle \sum _{k=1}^m
k\,\gamma_{k}  \,\tilde {u}_{{\tilde {x}}  }   =0
\end{aligned}
\end{equation}

\section{Solitary waves}
\label{Solitons}

\noindent Approximated solutions of the Burgers equation
(\ref{Burgers}) by means of the difference scheme (\ref{DRP})
strongly depend on the values of the time and space steps. For
specific values of $\tau$ and $h$, equation (\ref{EqEq}) can, for
instance, exhibit travelling wave solutions which can represent
great
disturbances of the searched solution.\\
\noindent We presently aim at determining the conditions, depending
on the values of the parameters $\tau$ and $h$, which give birth to
travelling wave solutions of (\ref{EqEq}).\\
\noindent Following Feng \cite{feng1} and our previous work
\cite{David}, in which travelling wave solutions of the
\textit{CBKDV} equation were exhibited as combinations of
bell-profile waves and kink-profile waves, we aim at determining
travelling wave solutions of (\ref{EqEq}) (see
\cite{li},\cite{whitham},
\cite{ablowitz}, \cite{dodd}, \cite{johnson}, \cite{ince}, \cite{zhang3}, \cite{birk}, \cite{Polyanin}).\\
\noindent Following \cite{feng1}, we assume that equation
(\ref{EqEq}) has  travelling wave solutions of the form
\begin{equation}\label{ChgtVar} \tilde{u}(\tilde{x}, \tilde{t}) =\tilde{u}(\xi), \quad \xi= \tilde{x}-v\,\tilde{t} \end{equation}
where $v$ is the wave velocity. Substituting (\ref{ChgtVar}) into
equation (\ref{ApproxAdim}) leads to:

\begin{equation}
\label{EqXi}  \widetilde{\mathcal {F}}
(\tilde{u},\tilde{u}^{(r)},(-v)^s\,\tilde{u}^{(s)}) = 0,
\end{equation}

\noindent Performing an integration of (\ref{EqXi}) with respect
to $\xi$ leads to an equation of the form:

\begin{equation}
\label{EqXiInt}  \widetilde{\mathcal {F}}_\xi^{\mathcal P}
(\tilde{u},\tilde{u}^{(r)},(-v)^s\,\tilde{u}^{(s)}) = C,
\end{equation}

\noindent where $C$ is an arbitrary integration constant, which
will be the starting point for the determination of solitary waves
solutions.\\

\noindent It is important to note that, contrary to other works,
the integration constant is not taken equal to zero, which would
lead to a loss of solutions.

\section{Travelling Solitary Waves}

\label{Solutions}

\subsection{Hyperbolic Ansatz}

\noindent The discussion in the preceding section provides us
useful information when we construct travelling solitary wave
solutions for equation (\ref{EqXi}). Based on these results, in
this section, a class of travelling wave solutions of the
equivalent equation (\ref{Eq}) is searched as a combination of
bell-profile waves and kink-profile waves of the form
\begin{equation}
\label{sol} \tilde{u}(\tilde{x}, \tilde{t}) = \sum_{i = 1}^n \left
(U_i\; \text{tanh}^i \left [\,C_i (\tilde{x}-v\,\tilde{t})\, \right
] + V_i \; \text{sech}^i \left [\,C_i(\tilde{x}-v\,\tilde{t}+x_0)\,
\right ] \right )+V_0
\end{equation}
where the $U_i's$, $V_i's$,  $C_i's$, $(i=1,\ \cdots,\ n)$, $V_0$
and $v$ are constants to be determined.
\\
\noindent In the following, $c$ is taken equal to 1.

\subsection{Theoretical analysis}

\noindent Substitution of (\ref{sol}) into equation
(\ref{EqXiInt}) leads to an equation of the form

\begin{equation}
\label{EqgenInt}
 \sum_{i,\, j, \,k} A_i\,\text{tanh}^i \big (C_i\,
\xi \big )\,\text{sech}^j \big (C_i \,\xi \big )\,\text{sinh}^k
\big (C_i\, \xi  \big )  =C
\end{equation}

\noindent the $A_i$ being real constants.\\
 \noindent The difficulty
for solving equation (\ref{EqgenInt}) lies in finding the values
of the constants $U_i$, $V_i$,  $C_i$, $V_0$ and $v$ by solving
the over-determined algebraic equations. Following \cite{feng1},
after balancing the higher-order derivative term and the leading
nonlinear term, we deduce $n=1$.\\

\noindent Then, following \cite{David} we replace
$\mbox{sech}({C_1} \,\xi)$ by $ \frac{2}{ e^{\,{C_1} \,\xi}+e^{\,-
{C_1} \,\xi } }$, $\mbox{sinh(}{C_1} \,\xi)$ by $ \frac {
e^{\,{C_1} \,\xi}-e^{\,- {C_1} \,\xi } }{2}$, $\mbox{tanh}({C_1}
\,\xi)$ by $\frac{e^{\, {C_1} \,\xi } -
      e^{\,- {C_1} \,\xi}}{e^{ \,{C_1} \,\xi   } + e^{\,-{C_1} \,\xi }}$, and multiply both sides by
      $(1+ e^{ 2\,\xi \, {C_1}})^2$, so that equation
(\ref{EqgenInt}) can be rewritten in the following form:
\begin {equation}
\label{Syst} \sum_{k=0}^{4} P_k ( U_1,\ V_1,\ C_1,\ v,\, V_0 )
\,e^{\,k \,C_1\, \xi} \; = \;0,
\end{equation}
where the $P_k$ $(k=0,\,...,\,4)$, are polynomials of $U_1$,
$V_1$, $C_1$, $V_0$  and $v$.

\noindent Depending wether (\ref{EqgenInt}) admits or no
consistent solutions,  spurious solitary waves solutions may, or
not, appear.

\subsection{Numerical scheme analysis}

\noindent Equation (\ref{EqXi}) is presently given by:

\begin{equation}
\label{EqXi1} \scriptsize  \begin{aligned} &  -v\,\tilde{u}'( \xi
) -\frac{v^2\,\sigma}{2 } \,\tilde {u}''(\xi)
  +\frac{ 2\,\sigma  }{\mu\,Re_h} \,\displaystyle \sum _{k=1}^m
k\,\gamma_{k}  \,\tilde {u}'(\xi) =0
\end{aligned}
\end{equation}
\noindent Performing an integration of (\ref{EqXi1}) with respect
to $\xi$ yields:

\begin{equation}
\label{EqXi2}  \scriptsize  \begin{aligned} & -v\,\tilde{u} ( \xi
) -\frac{v^2\,\sigma}{2 }\,\tilde {u}' (\xi)
  +\frac{ 2\,\sigma  }{\mu\,Re_h} \,\displaystyle \sum _{k=1}^m
k\,\gamma_{k}   \,\tilde {u} (\xi) =C
\end{aligned}
\end{equation}

\noindent i. e.:

\begin{equation}
\label{EqXi3}  \scriptsize    \left \lbrace \frac{ 2\,\sigma
}{\mu\,Re_h} \,\displaystyle \sum _{k=1}^m k\,\gamma_{k}  -v
\right \rbrace \,\tilde {u} (\xi) -\frac{v^2\,\sigma}{2 }\,\tilde
{u}' (\xi)=C
\end{equation}

\noindent where $C$ is an arbitrary integration constant.

\noindent Substitution of (\ref{sol}) for $n=1$ into equation
(\ref{EqXi3}) leads to:

\begin{equation}
\label{EqXi4}  \scriptsize    \left \lbrace \frac{ 2\,\sigma
}{\mu\,Re_h} \,\displaystyle \sum _{k=1}^m k\,\gamma_{k}  -v
\right \rbrace \,  \left \lbrace U_1\, \text{tanh}  \left [\,C_1
\,\xi\, \right ] + V_1 \, \text{sech} \left [\,C_1\,\xi\, \right ]
+V_0  \right \rbrace -\frac{v^2\,\sigma}{2 }\, \left \lbrace U_1\,
\text{sech}^2  \left [\,C_1\,\xi\, \right ] - V_1 \, \frac {
\text{sinh} \left [\,C_1\,\xi\, \right ] }{ \text{cosh}^2 \left
[\,C_1\,\xi\, \right ] } \right \rbrace=C
\end{equation}

\noindent i. e.:

\begin{equation}
\label{EqXi4}  \scriptsize    \left \lbrace \frac{ 2\,\sigma
}{\mu\,Re_h} \,\displaystyle \sum _{k=1}^mk\,\gamma_{k}  -v \right
\rbrace \,  \left \lbrace U_1\,
\frac{e^{C_1\,\xi}-e^{-C_1\,\xi}}{e^{C_1\,\xi}+e^{-C_1\,\xi}}\, +
  \frac{2\,V_1}{e^{C_1\,\xi}+e^{-C_1\,\xi}} +V_0  \right
\rbrace -\frac{v^2\,\sigma}{2 }\, \left \lbrace U_1\, \left (
\frac{2}{e^{C_1\,\xi}+e^{-C_1\,\xi}} \right )^2 - 2\,V_1 \,
\frac{e^{C_1\,\xi}-e^{-C_1\,\xi}}{\left (
{e^{\,C_1\,\xi}+e^{\,-\,C_1\,\xi}} \right )^2} \right \rbrace=C
\end{equation}

\noindent Multiplying both sides by $\left ( {1+e^{\,2\,C_1\,\xi}}
\right )^2$ yields:

\begin{equation}
\label{EqXi4}  \scriptsize
\begin{aligned} &  \left \lbrace \frac{ 2\,\sigma
}{\mu\,Re_h} \,\displaystyle \sum _{k=1}^mk\,\gamma_{k}  -v \right
\rbrace \,  \left \lbrace U_1\,
  \left ({e^{\,4\,C_1\,\xi}-1} \right ) \, +
   2\,V_1\,\left (e^{\,3\,C_1\,\xi}+e^{\,C_1\,\xi} \right ) +V_0 \,\left ( {1+e^{\,2\,C_1\,\xi}}
\right )^2 \right \rbrace\\
&-\frac{v^2\,C_1\,\sigma \,C_1}{2 }\, \left \lbrace 4\, U_1  -
2\,V_1 \,\left ({e^{\,3\,C_1\,\xi}-1}\right ) \right \rbrace=C
\end{aligned}
\end{equation}

\noindent which is a fourth-order equation in $e^{\,C_1\,\xi}$.
This equation being satisfied for any real value of $\xi$, one
therefore deduces that the coefficients of $e^{\,k\,C_1\,\xi}$,
$k=0,\, \ldots, \,4$ must be equal to zero, i.e.:

\begin{equation}
\label{EqXi4}  \scriptsize \left \lbrace
\begin{aligned}
 &   2\,\left \lbrace \frac{ 2\,\sigma
}{\mu\,Re_h} \,\displaystyle \sum _{k=1}^mk\,\gamma_{k}  -v \right
\rbrace \,  \left \lbrace -U_1 +V_0   \right
\rbrace-\frac{v^2\,C_1\,\sigma}{2 }\, \left \lbrace 4\, U_1  +
2\,V_1
  \right \rbrace=C\\
 &  \left
\lbrace \frac{ 2\,\sigma }{\mu\,Re_h} \,\displaystyle \sum
_{k=1}^mk\,\gamma_{k}  -v \right \rbrace \,
   2\,V_1     =0\\
&   2\,\left \lbrace \frac{ 2\,\sigma }{\mu\,Re_h} \,\displaystyle
\sum _{k=1}^mk\,\gamma_{k}  -v \right \rbrace \,
 V_0=0\\
&  2\, \left \lbrace \frac{ 2\,\sigma }{\mu\,Re_h} \,\displaystyle
\sum _{k=1}^m k\,\gamma_{k}  -v \right \rbrace \,
   V_1
 + v^2\,C_1\,\sigma \,V_1   =0\\
&  \left \lbrace \frac{ 2\,\sigma }{\mu\,Re_h} \,\displaystyle
\sum _{k=1}^mk\,\gamma_{k} -v \right \rbrace \, \left \lbrace U_1
     +V_0 \,
 \right \rbrace=0\\
\end{aligned} \right.
\end{equation}

\noindent $v=\frac{ 2\,\sigma }{\mu\,Re_h} \,\displaystyle \sum
_{k=1}^mk\,\gamma_{k} \, , \, V_1 \neq 0$ leads to the trivial nul
solution. Therefore, $V_1$ is necessarily equal to zero, which
implies:

\begin{equation}
\label{EqXi4}  \scriptsize \left \lbrace
\begin{aligned}
&v=\frac{ 2\,\sigma }{\mu\,Re_h} \,\displaystyle \sum
_{k=1}^mk\,\gamma_{k} \\
 &     U_1
   =-\frac{C}{2\,C_1\, v^2\,\sigma}\\
&  V_0 \in \R  \,\,\, , \,\,\,C_1 \in \R
\end{aligned} \right.
\end{equation}

\noindent All \textit{DRP} schemes admit thus kink-profile
travelling solitary waves solutions, given by:

\footnotesize

\begin{equation}
\label{sol} \tilde{u}(\tilde{x}, \tilde{t}) =  -\frac{C}{2\,C_1\,
\left ( \frac{ 2\,\sigma }{\mu\,Re_h} \,\displaystyle \sum
_{k=1}^mk\,\gamma_{k} \right ) ^2\,\sigma}\, \text{tanh} \left
[\,C_1 (\tilde{x}-v\,\tilde{t})\, \right ]  +V_0
\end{equation}

\normalsize

\section{Conclusions}

The analysis of the nonlinear equivalent differential equation for
finite-differenced \textit{DRP} schemes for the Burgers equation
has been carried out. We show that all \textit{DRP} schemes admit
spurious travelling solitary waves solutions, which make them, as
regards this point, structurally instable.

\addcontentsline{toc}{section}{\numberline{}References}


\begin{thebibliography}{1}

\normalsize{



\bibitem{Chaos1} David C., Sagaut P., {Spurious solitons and
structural stability of finite difference schemes for nonlinear
wave equations}, under press for \emph{Chaos, Solitons and
Fractals}.

\bibitem{burger1} Burgers J. M., Mathematical examples illustrating
relations occurring in the theory of turbulent fluid motion, {\em
Trans. Roy. Neth. Acad. Sci.} Amsterdam, 17 (1939) 1-53.


\bibitem{feng1} Feng Z. and Chen G., Solitary Wave Solutions of the Compound
Burgers-Korteweg-de Vries Equation, {\em Physica A}, 352 (2005)
419-435.


\bibitem{David}  David, Cl., Fernando, R., Feng Z., {A note on "general solitary wave solutions
of the Compound Burgers-Korteweg-de Vries Equation"}, Physica A:
Statistical and Theoretical Physics, 375 (1) (2007) 44-50.


\bibitem{Tam} Tam, Christopher K. W., Webb, Jay C., \emph{Dispersion-Relation-Preserving Finite Difference Schemes for
 Computational Acoustics}, Journal of Computational Physics, 107(2) (1993) 262-281.


\bibitem{Shokin} Shokin, Y. Liu, The method of differential approximation, Springer Verlag, Berlin (1983).


\bibitem{li} Li B., Chen Y. and Zhang H. Q., Explicit exact solutions for new general
two-dimensional KdV-type and two-dimensional KdV–Burgers-type
equations with nonlinear terms of any order, {\em J. Phys. A
(Math. Gen.)} 35 (2002) 8253–8265.

\bibitem{whitham} Whitham G. B., Linear and Nonlinear Waves,
Wiley-Interscience, New York, 1974.

\bibitem{ablowitz} Ablowitz M. J. and Segur H., Solitons and the Inverse Scattering
Transform, SIAM, Philadelphia, 1981.


\bibitem{dodd} Dodd R. K., Eilbeck J. C., Gibbon J.D. and Morris H. C., Solitons and Nonlinear Wave
Equations, London Academic Press, London, 1983.

\bibitem{johnson} Johnson R. S., A Modern Introduction to the Mathematical Theory of
Water Waves, Cambridge University Press, Cambridge, 1997.

\bibitem{ince} Ince E.L., Ordinary Differential Equations, Dover
Publications, New York, 1956.

\bibitem{zhang3} Zhang Z. F., Ding T.R., Huang W. Z. and Dong Z. X.,
Qualitative Analysis of Nonlinear Differential Equations, Science
Press, Beijing, 1997.

\bibitem{birk} Birkhoff G. and Rota G. C., Ordinary Differential
Equations, Wiley, New York, 1989.


\bibitem{Polyanin} Polyanin A. D. and Zaitsev V. F., Handbook of Nonlinear Partial Differential
Equations, Chapman and Hall/CRC, 2004. }


\end{thebibliography}
\end{document}